# A HOMING PROBLEM FOR DIFFUSION PROCESSES WITH CONTROL-DEPENDENT VARIANCE[1]

BY MARIO LEFEBVRE

*École Polytechnique de Montréal*

Controlled one-dimensional diffusion processes, with infinitesimal variance (instead of the infinitesimal mean) depending on the control variable, are considered in an interval located on the positive half-line. The process is controlled until it reaches either end of the interval. The aim is to minimize the expected value of a cost criterion with quadratic control costs on the way and a final cost equal to zero (resp. a large constant) if the process exits the interval through its left (resp. right) end point. Explicit expressions are obtained both for the optimal value of the control variable and the value function when the infinitesimal parameters of the processes are proportional to a power of the state variable.

**1. Introduction.** Let $x(t)$ be a one-dimensional controlled diffusion process defined by the stochastic differential equation

$$dx(t) = a[x(t), t] \, dt + b[x(t), t] u(t) \, dt + \{v[x(t), t]\}^{1/2} \, dW(t),$$

where $v(\cdot, \cdot)$ is positive and $W(t)$ is a standard Brownian motion, and define

$$\tau(x) := \inf\{t \geq 0 : x(t) \in \{d_1, d_2\} \mid x(0) = x \in [d_1, d_2]\}.$$

Using a result due to Whittle [(1982), page 289], we can show that if the uncontrolled process $y(t)$ that corresponds to $x(t)$ is certain to leave the interval $[d_1, d_2]$, then the value of the control $u(t)$ that minimizes the expected value of the cost criterion

$$J(x) = \int_0^{\tau(x)} \tfrac{1}{2} q[x(t), t] u^2(t) \, dt + K(x(\tau), \tau),$$

Received September 2001; revised May 2003.

[1]Supported by the Natural Sciences and Engineering Research Council of Canada.

where $q(\cdot,\cdot)$ is positive and $K(\cdot,\cdot)$ is a general termination cost function, can be obtained from the mathematical expectation

$$E[e^{-K(y(\tau),\tau)/\alpha}|y(0)=x],$$

in which $\alpha$ is a positive parameter such that the relation

$$v = \alpha b^2/q$$

holds. Actually, the mathematical expectation above is equal to $e^{-F(x)/\alpha}$, where $F(x)$ is the value function defined by

(1) $$F(x) = \inf_{u(t),0 \leq t \leq \tau(x)} E[J(x)].$$

Whittle has termed this type of problem $LQG$ homing. It can be generalized [see Lefebvre (1989, 1997)] in particular by using a risk-sensitive cost criterion instead of $J(x)$ [see Kuhn (1985) and Whittle (1990), page 222].

Next, the author [Lefebvre (2001)] modified the problem set up by Whittle by considering the controlled process defined by

$$dx(t) = b[x(t)]u(t)\,dt + \{v[x(t)]|u(t)|\}^{1/2}\,dW(t)$$

and the cost criterion

$$J(x) = \int_0^{\tau(x)} \{\tfrac{1}{2}q[x(t)]u^2(t) + \lambda\}\,dt,$$

where $q(\cdot)$ is positive, $\lambda$ is a positive parameter and

$$\tau(x) := \inf\{t \geq 0 : |x(t)| = d \mid x(0) = x \in [-d,d]\}.$$

He found, under some symmetry assumptions, that the control that minimizes the expected value of $J(x)$ is given by

$$u^* = \left(\frac{2\lambda}{q(x)}\right)^{1/2} \qquad \text{when } 0 \leq x \leq d.$$

He also gave a probabilistic interpretation to the value function $F(x)$ and, finally, he computed explicitly this function in the most important cases, for instance, the cases when $x(t)$ with $u(t) \equiv 1$ is a Wiener process or a geometric Brownian motion.

In the present paper, we assume that the controlled stochastic process $X(t)$ obeys the stochastic differential equation

(2) $$dX(t) = f[X(t)]\,dt + \{v[X(t)]|u(t)|\}^{1/2}\,dW(t),$$

where $f(x)$ and $v(x)$ are positive functions for $x > d_1$, and we let

$$\tau(x) := \inf\{t \geq 0 : X(t) \in \{d_1,d_2\} \mid X(0) = x \in [d_1,d_2]\}$$

with $d_1 \geq 0$. Therefore, it is the infinitesimal variance of the controlled process $X(t)$ that is control-dependent, rather than its infinitesimal mean, as in Whittle [(1982), page 289]. The cost criterion is

$$J(x) = \int_0^{\tau(x)} \{\tfrac{1}{2} q[X(t)] u^2(t) + \lambda\} \, dt + K\{X[\tau(x)]\}. \tag{3}$$

In $J(x)$, the parameter $\lambda$ can now take any real value and $K(\cdot)$ is defined by

$$K\{X[\tau(x)]\} = \begin{cases} 0, & \text{if } X[\tau(x)] = d_1, \\ K_0, & \text{if } X[\tau(x)] = d_2, \end{cases}$$

where $K_0$ is a (large enough) positive constant. Thus, if $\lambda$ is positive, then the aim is to make the controlled process $X(t)$ leave the interval $(d_1, d_2)$ as soon as possible and through its left end point $d_1$, whereas when $\lambda$ is negative there is a reward given for survival in $(d_1, d_2)$. If $\lambda = 0$, time spent in the continuation region $(d_1, d_2)$ is neither directly rewarded nor penalized; however, because the function $q(\cdot)$ is assumed to be strictly positive, the sooner $X(t)$ exits $(d_1, d_2)$, the better. In all cases, the quadratic control costs must of course be taken into account.

In the next section, the optimal value of the control variable will be computed. The case when the functions $f$, $v$ and $q$ are proportional to $X^n(t)$, where $n \in \{0, 1, 2\}$, will be treated thoroughly. Particular examples, including the case when $X(t)$ with $u(t) \equiv 1$ is a Wiener process, will be presented in Section 3. Finally, a few concluding remarks will be made in Section 4.

**2. Computation of the optimal control.** Let $F(x)$ be the value function defined in (1). Assuming that it exists and is twice differentiable, we can easily show that it satisfies the dynamic programming equation

$$0 = \inf_u \{\tfrac{1}{2} q(x) u^2 + \lambda + f(x) F'(x) + \tfrac{1}{2} v(x) |u| F''(x)\} \tag{4}$$

for $d_1 < x < d_2$, where $u := u(0)$. Moreover, the boundary conditions are

$$F(d_1) = 0 \quad \text{and} \quad F(d_2) = K_0. \tag{5}$$

Now, since $u(t)$ only appears in absolute value in the stochastic differential equation (2) and squared in the cost criterion (3), the sign of $u(t)$ is actually irrelevant. Hence, we can assume without loss of generality that $u$ is nonnegative and it then follows at once that the optimal control $u^*$ is given by

$$u^* = -\frac{v(x)}{2q(x)} F''(x) \ (\geq 0). \tag{6}$$

REMARKS. 1. We will have to check below that $F''(x)$ is indeed less than or equal to zero. Actually, we cannot have $F''(x) = 0$ and satisfy both boundary conditions in (5). So, $F''(x)$ should in fact be strictly negative.
2. We have assumed above that the function $f(x)$ is positive if $x > d_1$. Notice that if $f(x)$ is negative for $x \in [d_1, d_2]$ and $\lambda = 0$, then the optimal control is trivially given by $u^* \equiv 0$. Indeed, we then obtain that $F(x) = 0$ (for $d_1 \leq x < d_2$), which is clearly the smallest value that $F(x)$ can take. However, when we choose a parameter $\lambda$ different from zero, the case $f(x) \neq 0$ could be considered.
3. If the parameter $\lambda = 0$, then the function $F(x)$ takes its values in the interval $[0, K_0]$.
4. In some cases, the origin is an inaccessible boundary for the uncontrolled process $X_1(t)$ obtained by setting $u(t) \equiv 1$; that is, the origin cannot be reached in finite time. This is true, in particular, when $X_1(t)$ is a geometric Brownian motion defined by the stochastic differential equation

$$dX_1(t) = X_1(t)\,dt + |X_1(t)|\,dW(t).$$

Therefore, it is natural to choose $d_1$ strictly positive in such a case. However, in other cases $d_1$ can be chosen equal to zero in a very legitimate way.

Next, substituting the optimal control $u^*$ into the dynamic programming equation (4), we find that the value function $F(x)$ satisfies the second-order nonlinear ordinary differential equation

(7) $$0 = \lambda + f(x)F'(x) - \frac{v^2(x)}{8q(x)}[F''(x)]^2.$$

Summing up, we may state the following proposition.

PROPOSITION 2.1. *If the value function $F(x)$ exists and is twice differentiable, then the optimal control $u^*$ is given by the formula* (6). *Moreover, the function $F(x)$ can be obtained by solving the ordinary differential equation* (7), *subject to the boundary conditions* (5).

REMARK. Note that for the problem set up above to make sense, we must have $\lambda + f(x)F'(x) \geq 0$ [see (7)].

Next, we will solve explicitly the nonlinear ordinary differential equation (7) in two particular cases.

CASE 1. Assume first that $f(x) \equiv f_0$, a positive constant, and that

(8) $$\frac{v^2(x)}{8q(x)} := h(x) = h_0 x^n,$$

where $h_0$ is also a positive constant and $n \in \{-2, -1, \ldots, 4\}$.

REMARK. The most important cases for applications are the ones when $f(x) = f_0 x^k$, $v(x) = v_0 x^j$ and $q(x) = q_0 x^l$ with $j, k, l \in \{0, 1, 2\}$. Actually, we could also include the case when $f(x) = f_0/x$. Indeed, $X(t)$ with $u(t) \equiv 1$ could then be a Bessel process [if $v(x) \equiv 1$; see Karlin and Taylor (1981), pages 175 and 176, for instance].

When the formula (8) holds [and $f(x) \equiv f_0$], we may rewrite the ordinary differential equation (7) as

(9) $$h(x)[F''(x)]^2 = \lambda + f_0 F'(x).$$

Differentiating this differential equation, we obtain that

$$h'(x)[F''(x)]^2 + 2h(x)F''(x)F'''(x) = f_0 F''(x).$$

Hence, because we must have $F''(x) < 0$ (see above), we may write that

(10) $$2h(x)G'(x) + h'(x)G(x) = f_0,$$

where $G(x) := F''(x)$. The general solution of (10) is given by

$$G(x) = h^{-1/2}(x)\left\{f_0 \int \frac{1}{2h^{1/2}(x)} \, dx + c\right\},$$

where $c$ is a constant. It is now easy to obtain an explicit expression in the case when $h(x) = h_0 x^n$. We find that

$$G(x) = \frac{c}{\sqrt{h_0}} x^{-n/2} + \frac{f_0}{h_0} \frac{x^{1-n}}{2-n} \qquad \text{if } n \neq 2$$

(11) and

$$G(x) = \frac{c}{\sqrt{h_0}} \frac{1}{x} + \frac{f_0}{h_0} \frac{\ln x}{2x} \qquad \text{if } n = 2.$$

Integrating the function $G(x)$ twice, we obtain that the value function $F(x)$ is given by

(12) $$\begin{cases} \dfrac{2c}{\sqrt{h_0}} \dfrac{x^{2-(n/2)}}{(2-n)(2-n/2)} + \dfrac{f_0}{h_0} \dfrac{x^{3-n}}{(2-n)^2(3-n)} + c_1 x + c_0 & \\ \hspace{6cm} \text{if } n \neq 2, 3, 4, \\[6pt] \dfrac{c}{\sqrt{h_0}} x(\ln x - 1) + \dfrac{f_0}{4h_0} x(\ln^2 x - 2\ln x + 2) + c_1 x + c_0 & \\ \hspace{6cm} \text{if } n = 2, \\[6pt] -\dfrac{c}{\sqrt{h_0}} 4 x^{1/2} + \dfrac{f_0}{h_0} \ln x + c_1 x + c_0 & \text{if } n = 3, \\[6pt] -\dfrac{c}{\sqrt{h_0}} \ln x - \dfrac{f_0}{4h_0} \dfrac{1}{x} + c_1 x + c_0 & \text{if } n = 4, \end{cases}$$

where $c_1$ and $c_0$ are constants.

Finally, the constants $c$, $c_1$ and $c_0$ are uniquely determined from the boundary conditions (5) and the equation (9). Actually, we find that $c^2 = \lambda + f_0 c_1$. An example will be presented in the next section.

CASE 2. If $\lambda = 0$ and
$$\frac{v^2(x)}{8q(x)f(x)} := g(x) = g_0 x^m,$$

where $g_0$ is a positive constant and $m \in \{-4, -3, \ldots, 4\}$, we have

(13) $$g(x)[F''(x)]^2 = F'(x).$$

Proceeding as above, we find that the function $F(x)$ is given by (12), with $f_0 = 1$, $n$ replaced by $m$ and $h_0$ by $g_0$, and that the constants $c$, $c_1$ and $c_0$ are now uniquely determined from the boundary conditions (5) and the ordinary differential equation (13). Here, we find that $c^2 = c_1$. As for Case 1, an example will be provided in the next section.

REMARKS.
1. As mentioned above, we must also check that the condition $F''(x) < 0$ is satisfied.
2. As will be seen in the examples presented in Section 3, there is also a restriction on the constant $K_0$ in the definition of the function $K(\cdot)$.

**3. Examples.** (a) First, we consider the particular case when $f(x) \equiv f_0$, $v(x) \equiv v_0$ and $q(x) \equiv q_0$, where $f_0$, $v_0$ and $q_0$ are all positive constants. Then, the controlled process $X(t)$ with $u(t) \equiv 1$ is a Wiener process with infinitesimal parameters $f_0$ and $v_0$. We have

$$h(x) = h_0 = \frac{v_0^2}{8q_0},$$

so that $n = 0$ in (8). It follows that [see (12)]

$$F(x) = \frac{2f_0 q_0}{3v_0^2}x^3 + \frac{c\sqrt{2q_0}}{v_0}x^2 + c_1 x + c_0.$$

We find that the ordinary differential equation (9) is satisfied if and only if we take $c^2 = \lambda + f_0 c_1$, as noticed above. It follows that

$$F(x) = \frac{2f_0 q_0}{3v_0^2}x^3 + \frac{c\sqrt{2q_0}}{v_0}x^2 + \frac{c^2 - \lambda}{f_0}x + c_0.$$

Next, the boundary condition $F(d_1) = 0$ implies that

$$F(x) = \frac{2f_0 q_0}{3v_0^2}(x^3 - d_1^3) + \frac{c\sqrt{2q_0}}{v_0}(x^2 - d_1^2) + \frac{c^2 - \lambda}{f_0}(x - d_1),$$

and finally, $F(d_2) = K_0$ yields that

$$K_0 = \frac{c^2(d_2 - d_1)}{f_0} + \frac{c\sqrt{2q_0}(d_2^2 - d_1^2)}{v_0} + \frac{2f_0 q_0 (d_2^3 - d_1^3)}{3v_0^2} - \frac{\lambda(d_2 - d_1)}{f_0}.$$

Hence, we have

$$c = \frac{f_0}{2(d_2 - d_1)}$$
$$\times \left[ -\frac{\sqrt{2q_0}(d_2^2 - d_1^2)}{v_0} \right.$$
$$\left. \pm \left\{ \frac{2q_0(d_2^2 - d_1^2)^2}{v_0^2} \right.\right.$$
$$\left.\left. - 4\frac{d_2 - d_1}{f_0}\left(\frac{2f_0 q_0(d_2^3 - d_1^3)}{3v_0^2} - \frac{\lambda(d_2 - d_1)}{f_0} - K_0\right) \right\}^{1/2} \right].$$

To simplify further, we take $f_0 = v_0 = 1$ and $q_0 = 1/2$. We get that

(14)
$$c = -\frac{(d_2 + d_1)}{2}$$
$$\pm \frac{\{(d_2^2 - d_1^2)^2 - 4(d_2 - d_1)(1/3(d_2^3 - d_1^3) - \lambda(d_2 - d_1) - K_0)\}^{1/2}}{2(d_2 - d_1)}.$$

We then deduce that the constant $K_0$ must satisfy the inequality

(15) $$K_0 \geq \lambda(d_2 - d_1) + \tfrac{1}{12}(d_2 - d_1)^3.$$

Notice that the larger the parameter $\lambda$ is, the larger $K_0$ must be. Conversely, if $\lambda$ is (negative and) small enough, any $K_0 > 0$ is admissible.

To determine whether we must choose the "+" or "−" sign in (14), we will use the fact that we must have $F''(x) < 0$; that is, with $f_0 = v_0 = 1$ and $q_0 = \tfrac{1}{2}$,

$$F''(x) \equiv G(x) = 2x + 2c < 0 \qquad \text{for } d_1 \leq x \leq d_2.$$

This implies that the constant $c$ must be smaller than or equal to $-d_2$. Using (14), we obtain that

(16)
$$\pm \frac{\{(d_2^2 - d_1^2)^2 - 4(d_2 - d_1)(1/3(d_2^3 - d_1^3) - \lambda(d_2 - d_1) - K_0)\}^{1/2}}{2(d_2 - d_1)}$$
$$< \frac{d_1 - d_2}{2} \ (< 0).$$

Thus, we must choose the "−" sign. Then, we find that (16) implies that

(17)
$$(d_2^2 - d_1^2)^2 - 4(d_2 - d_1)(\tfrac{1}{3}(d_2^3 - d_1^3) - \lambda(d_2 - d_1) - K_0) > (d_2 - d_1)^4$$
$$\iff K_0 > \lambda(d_2 - d_1) + \tfrac{1}{3}(d_2 - d_1)^3.$$

Since this last condition is more restrictive than the one in (15), this is a condition that must be fulfilled.

Summing up, when we choose $f(x) \equiv 1$, $v(x) \equiv 1$ and $q(x) \equiv \frac{1}{2}$, the optimal control is given by

$$u^* = -2(x+c),$$

where the constant $c$ is defined in (14), in which the "$-$" sign is chosen and the constant $K_0$ satisfies the condition (17). Furthermore, the value function is

$$F(x) = \tfrac{1}{3}(x^3 - d_1^3) + c(x^2 - d_1^2) + (c^2 - \lambda)(x - d_1) \qquad \text{for } d_1 \le x \le d_2.$$

(b) Next, we consider the controlled stochastic process $X(t)$ defined by the stochastic differential equation

$$dX(t) = X(t)\, dt + \{X^2(t)|u(t)|\}^{1/2}\, dW(t)$$

and we look for the control $u^*$ that minimizes the expected value of the cost criterion

$$J(x) = \int_0^{\tau(x)} \tfrac{1}{2} X(t) u(t)\, dt + K\{X[\tau(x)]\}.$$

That is, we take $f[X(t)] = X(t)$, $v[X(t)] = X^2(t)$ and $q[X(t)] = X(t)$. Moreover, we set $\lambda = 0$. In the case when $u(t) \equiv 1$, $X(t)$ is a geometric Brownian motion. Notice that

$$\frac{v^2(x)}{8q(x)f(x)} = \frac{x^2}{8} \equiv g(x).$$

It follows that the value function is given by [see (12) with $n\ (=m) = 2$, $f_0 = 1$ and $h_0\ (=g_0) = \frac{1}{8}$]

$$F(x) = 2\sqrt{2}cx(\ln x - 1) + 2x(\ln^2 x - 2\ln x + 2) + c_1 x + c_0$$

and

$$G(x) \equiv F''(x) = 2\sqrt{2}c\frac{1}{x} + 4\frac{\ln x}{x}.$$

Furthermore, the optimal control $u^*$ is

(18) $$u^* = -\sqrt{2}c - 2\ln x.$$

Next, (13) yields that

$$c^2 + 2\ln^2 x + 2\sqrt{2}c\ln x = 2\sqrt{2}c\ln x + 2\ln^2 x + c_1 \implies c^2 = c_1,$$

which is in fact true for all problems in Case 2. Hence, we have

$$F(x) = 2\sqrt{2}cx(\ln x - 1) + 2x\ln x(\ln x - 2) + 4x + c^2 x + c_0.$$

The boundary conditions $F(d_1) = 0$ and $F(d_2) = K_0$ imply that

$$F(x) = 2\sqrt{2}c[(x \ln x - d_1 \ln d_1) - (x - d_1)] + 2(x \ln^2 x - d_1 \ln^2 d_1)$$
(19)
$$- 4(x \ln x - d_1 \ln d_1) + 4(x - d_1) + c^2(x - d_1),$$

where the constant $c$ is such that

$$K_0 = c^2(d_2 - d_1) + 2\sqrt{2}c[(d_2 \ln d_2 - d_1 \ln d_1) - (d_2 - d_1)]$$
$$+ 2d_2 \ln^2 d_2 - 2d_1 \ln^2 d_1 - 4d_2 \ln d_2 + 4d_1 \ln d_1 + 4(d_2 - d_1).$$

For simplicity, let $d_1 = 1$ and $d_2 = d$. Then, we have

$$0 = (d-1)c^2 + 2\sqrt{2}c[d \ln d - (d-1)] + 2d \ln^2 d - 4d \ln d + 4(d-1) - K_0,$$

so that

(20) $$c = \frac{-2\sqrt{2}(d \ln d - d + 1) \pm [8d \ln^2 d - 8(d-1)^2 + 4(d-1)K_0]^{1/2}}{2(d-1)}.$$

Thus, a first restriction on the constant $K_0$ is that

$$K_0 \geq 2(d-1) - 2\frac{d}{d-1} \ln^2 d,$$

where $d > 1$.

Finally, we know that we must also have

$$F''(x) = 2\sqrt{2}c\frac{1}{x} + 4\frac{\ln x}{x} < 0 \qquad \forall x \in [1, d],$$

which implies that

$$c < -\sqrt{2} \ln d.$$

We then deduce from (20) that we must again choose the "$-$" sign and that the constant $K_0$ must satisfy the inequality

$$K_0 > 4(d-1) - 4 \ln d - 2 \ln^2 d.$$

Since

$$4(d-1) - 4\ln d - 2\ln^2 d > 2(d-1) - 2\frac{d}{d-1}\ln^2 d \qquad \forall d > 1,$$

we must impose the constraint

(21) $$K_0 > 4(d - 1 - \ln d - \tfrac{1}{2}\ln^2 d).$$

In summary, when we choose $f[X(t)] = X(t)$, $v[X(t)] = X^2(t)$, $q[X(t)] = X(t)$, $\lambda = 0$, $d_1 = 1$ and $d_2 = d$ ($> 1$), the optimal control is given by (18) with the constant $c$ defined in (20), in which the "$-$" sign is chosen and the constant $K_0$ is such that (21) holds. Moreover, the value function is [see (19) with $d_1 = 1$ and $d_2 = d$]

$$F(x) = 2x \ln x[\ln x + \sqrt{2}c - 2] + (x - 1)[c^2 - 2\sqrt{2}c + 4].$$

**4. Concluding remarks.** In this article, the problem of optimally controlling a certain class of one-dimensional diffusion processes was set up and solved exactly. Contrary to the classic formulation, the diffusion processes in question had control-dependent infinitesimal variances rather than infinitesimal means. The objective, when the parameter $\lambda$ in the cost criterion $J(x)$ defined in (3) is nonnegative, was to incite the controlled process to exit, as soon as possible, the interval $[d_1, d_2]$ at $d_1$. When $\lambda$ is negative, a reward is given for survival in the interval $(d_1, d_2)$.

In Section 3, two particular cases were treated extensively. The first example involved a controlled Wiener process, whereas in the second example a geometric Brownian motion was optimally controlled. Many other important cases could be considered. For instance, particular Bessel processes could have been used. It would also be interesting to take $f[X(t)] = -\alpha X(t)$, with $\alpha$ a positive constant, $v[X(t)] \equiv v_0 > 0$ and $\lambda < 0$. Then, the controlled process $X(t)$ with $u(t) \equiv 1$ is an Ornstein–Uhlenbeck process. Since $f[X(t)]$ is negative when $X(t)$ is in the interval $[d_1, d_2]$, we have $J(x) = 0$ if $\lambda = 0$ and if the optimizer chooses $u(t) \equiv 0$. However, if $\lambda$ is negative and small enough (i.e., large enough in absolute value), the optimizer can receive a reward overall if the controlled process remains between $d_1$ and $d_2$ for a long enough time. Therefore, the optimal control $u^*$ will not always be identical to zero if $\lambda$ is small and $K_0$ is not too large. Actually, the situation is similar even when $f[X(t)] > 0$ in the interval $[d_1, d_2]$. Indeed, if $K_0$ is not large, the optimizer is better off to let $X(t)$ hit $d_2$ before $d_1$ rather than to use a lot of control to make $X(t)$ hit $d_1$ first; hence the constraints that we must impose on the constant $K_0$, as we have seen in the examples presented in Section 3.

Finally, possible extensions of the work presented in this article are the following: first, a two-dimensional version of the optimal control problem could be considered. Also, we could use a risk-sensitive cost criterion rather than $J(x)$.

## REFERENCES

Karlin, S. and Taylor, H. (1981). *A Second Course in Stochastic Processes*. Academic Press, New York. MR611513

Kuhn, J. (1985). The risk-sensitive homing problem. *J. Appl. Probab.* **22** 167–172. MR808859

Lefebvre, M. (1989). An extension of a stochastic control theorem due to Whittle. *IEEE Trans. Automat. Control* **34** 567–568. MR991909

Lefebvre, M. (1997). Reducing a nonlinear dynamic programming equation to a Kolmogorov equation. *Optimization* **42** 125–137. MR1488096

Lefebvre, M. (2001). A different class of homing problems. *Systems Control Lett.* **42** 347–352. MR2006862

Whittle, P. (1982). *Optimization Over Time* **1**. Wiley, Chichester. MR670949

Whittle, P. (1990). *Risk-Sensitive Optimal Control*. Wiley, Chichester. MR1093001

Département de Mathématiques
 et de Génie Industriel
École Polytechnique de Montréal
C.P. 6079, Succursale Centre-ville
Montréal, Québec
Canada H3C 3A7
e-mail: mlefebvre@polymtl.ca